\documentclass{amsproc}

\usepackage{amsfonts,amsmath,amssymb,amsthm}
\usepackage{mathrsfs}
\usepackage{color}

\usepackage[colorlinks,backref]{hyperref}

\newtheorem{theorem}{Theorem}[section]
\newtheorem{lemma}[theorem]{Lemma}
\newtheorem{prp}[theorem]{Proposition}

\theoremstyle{definition}

\newtheorem{example}[theorem]{Example}
\newtheorem{cor}[theorem]{Corollary}

\theoremstyle{remark}

\numberwithin{equation}{section}

\DeclareMathOperator{\Ric}{Ric} \DeclareMathOperator{\grad}{grad}

\copyrightinfo{}% copyright year
{}% copyright holder

\begin{document}

\title{Hessian Tensor and Standard Static Space-times}

\author{Fernando Dobarro}
%\address{Department of Mathematics and Informatics, University of Trieste, Trieste, Italy}
\address{
Dipartimento di Matematica e Informatica, Universit\`{a} degli
Studi di Trieste, Via Valerio 12/B, I-34127 Trieste, Italy}
%\curraddr{}
\email{dobarro@dmi.units.it}
%\thanks{The first author was supported in part by NSF Grant \#000000.}
\author{B\"{u}lent \"{U}nal}
\address{Department of Mathematics, Bilkent University,
         Bilkent, 06800 Ankara, Turkey}
\email{bulentunal@mail.com}
%\thanks{Support information for the second author.}
\subjclass[2000]{53C21, 53C50, 53C80}
\date{February 26, 2008}
%\dedicatory{This paper is dedicated to our advisors.}
\keywords{Warped products, Hessian tensor, Killing vector fields,
energy conditions,
%conformal hyperbolicity, conjugate points, time-like diameter,
standard static space-times.}

%\begin{abstract}
%We survey applications of Hessian tensor type equations on
%standard static space-times in terms of characterization Killing
%vector fields, energy conditions, conformal hyperbolicity and
%conjugate points along causal geodesics.
%\end{abstract}

\begin{abstract}
In this brief survey, we will remark the interaction among the
Hessian tensor on a semi-Riemannian manifold and some of the
several questions in Lorentzian (and also in semi-Riemannian)
geometry where this $2-$covariant tensor is involved. In
particular, we deal with the characterization of Killing vector
fields and the study of a set of consequences of energy conditions
in the framework of standard static space-times.
\end{abstract}

\maketitle

\section{Introduction}

The two central concepts in this note will be the Hessian tensor
on a semi-Riemannian manifold and the warped product of
semi-Riemannnian manifolds (specially standard static
space-times). There are many arguments in mathematical-physics
where these concepts interact.

We briefly recall some basic definitions. Let $(F,g_F)$ be a
semi-Riemannian manifold and $\varphi \in C^\infty(F)$ be a smooth
function on $F.$ Then the \textit{Hessian of $\varphi$} is the
$(0,2)-$tensor defined by \begin{equation}\label{eq:def-hessian}
   {\rm H}_F^\varphi(X,Y) =
   g_F(\nabla_X^F \grad_F \varphi,Y)= \nabla^F \nabla^F \varphi (X,Y) ,
\end{equation}
for any vector fields $X,Y \in \mathfrak{X} (F)$
where $\nabla^F$ is the Levi-Civita connection and $\grad_F$ is
the $g_F-$gradient operator. The $g_F-$trace of ${\rm
H}_F^\varphi$ is the Laplace-Beltrami operator,
$\Delta_{F}\varphi$. Notice that $\Delta_{F}$ is
%uniformly
elliptic if $(F,g_F)$ is Riemannian.

Let $(B,g_B)$ and $(F,g_F)$ be pseudo-Riemannian manifolds and
also let $b \colon B \to (0,\infty)$ be a smooth function. Then
the (singly) \textit{warped product}, $B \times_b F$ is the
product manifold $B \times F$ furnished with the metric tensor
$g=g_B \oplus b^{2}g_F$ defined by
\begin{equation}\label{eq:wp-metric}
g=\pi^{\ast}(g_B) \oplus (b \circ \pi)^2 \sigma^{\ast}(g_F),
\end{equation}
where $\pi \colon B \times F \to B$ and $\sigma \colon B \times F
\to F$ are the usual projection maps and ${}^\ast$~denotes the
pull-back operator on tensors. Here, the function $b$ is called
the warping function. Warped product manifolds were introduced in
general relativity as a method to find general solutions to
Einstein's field equations \cite{B,BEE,ON}. Two important examples
include generalized Robertson-Walker space-times and standard
static space-times (a generalization of the Einstein static
universe). Precisely a \textit{standard static space-time} is a
Lorentzian warped product where the warping function is defined on
a Riemannian manifold called the base and acting on the negative
definite metric on an open interval of real numbers, called the
fiber. More precisely, a standard static space-time, denoted by
$I_f \times F$, is a Lorentzian warped product furnished with the
metric $g=-f^2{\rm d}t^2 \oplus g_F,$ where $(F,g_F)$ is a
Riemannian manifold, $f \colon F \to (0,\infty)$ is smooth and
$I=(t_1,t_2)$ with $-\infty \leq t_1 < t_2 \leq \infty$. In
\cite{ON}, it was shown that any static space-time is locally
isometric to a standard static space-time.

There are many subjects in semi-Riemannian geometry and physics
where all these ingredients interact and play a central role. For
instance in the study of concircular scalar fields
\cite{Obata1962,T65}; in recent studies of Hessian manifolds
\cite{Shima2007}; in several questions of curvature of warped
products and the construction of Einstein manifolds
\cite{B,DU04a,DU07a,DU07b} and in the characterization of Killing
vector fields on Robertson Walker space-times \cite{S99}, among
many others. We will concentrate our attention to the study of
Killing vector fields and energy conditions on standard static
space-times, where the ingredients mentioned above are involved
(see \cite{DU08a} and \cite{DU08b}). The references mentioned above
are mere indications in which the reader can find more specific
references for each argument as well as links to alternative
approaches and much more.

\section{Preliminaries and Notation}

Throughout the paper $I$ will denote an open real interval
$I=(t_1,t_2)$, where $-\infty \leq t_1 < t_2 \leq \infty.$
Moreover, $(F,g_F)$ will be a connected Riemannian manifold
without boundary with $\dim F=s$. Finally, on an arbitrary
differentiable manifold $N$, $C^{\infty}_{>0}(N)$ denotes the set
of all strictly positive $C^{\infty}$ functions defined on $N$ and
$\mathfrak X(N)$ will denote the $C^{\infty}(N)-$module of smooth
vector fields on $N$.

Suppose that $U \in \mathfrak X (I)$ and $V,W \in \mathfrak
X (F)$. If ${\rm Ric}$ and ${\rm Ric}_F$ denote the Ricci tensors
of $I _f\times F$ and $(F,g_F)$, respectively, then
\begin{equation} \label{eq:ricci-sss}
{\rm Ric} \left(U+V, U+W \right)={\rm Ric}_F(V,W)+ f \Delta_F f \,
{\rm d}t^2(U,U) - \frac{1}{f}{\rm H}^f_F(V,W).
\end{equation}

If $\tau$ and $\tau_F$ denote the scalar curvatures of $I _f\times
F$ and $F$, respectively, then
\begin{equation} \label{eq:sc-sss} \tau = \tau_F -2
\frac{1}{f}\Delta_F f.
\end{equation}

From now on, for any given $f \in C^{\infty}_{>0}(F)$, $Q_F^f$ will
denote the $2$ covariant tensor
\begin{equation}\label{eq: quadratic form 1}
    Q_F^f:=\Delta_F f \, g_F- H_F^f.
\end{equation}

%The quadratic forms associated to $\Ric_F$ and $Q^f_F$ will be
%denoted by $\mathcal{R}ic_F$ and $\mathcal{Q}^f_F$, respectively.

$\mathcal{R}ic_F$(respectively, $\mathcal{Q}^f_F$) denotes the quadratic
form associated to $\Ric_F$ (respectively, $Q^f_F$).

Notice that \eqref{eq:ricci-sss} and \eqref{eq: quadratic form 1}
imply that for any $U \in \mathfrak X (I)$ and $V,W \in \mathfrak
X (F)$ is
\begin{equation}\label{eq:Qricci-sss}
    {\rm Ric} (U+V,U+W)= {\rm Ric}_F(V,W) + \frac{1}{f}
    Q_F^f(V,W)- g(U+V,U+W) \frac{1}{f}\Delta_F f.
\end{equation}

\section{Killing Vector Fields}
To begin with, we recall the concepts of Killing and conformal-Killing
vector fields on pseudo-Riemannian manifolds. Let $(N, g_N)$ be a
pseudo-Riemannian manifold and $X \in \mathfrak X(N)$. Then
\begin{itemize}
\item $X$ is said to be Killing if ${\rm L}_X g_N = 0$,

\item $X$ is said to be conformal-Killing if $\exists \sigma \in
C^{\infty}(N)$ such that ${\rm L}_X g_N = 2 \sigma  g_N$,
\end{itemize}
where ${\rm L}_X$ denotes the Lie derivative with respect to $X$.
Moreover, for any $Y$ and $Z$ in $\mathfrak X(N),$ we have the
following identity (see \cite[p.250 and p.61]{ON})
\begin{equation}\label{eq:Lie deriv 1}
    {\rm L}_X g_N(Y,Z) = g_N(\nabla_Y X,Z) + g_N(Y,\nabla_Z X).
\end{equation}

\noindent Notice that, any vector field on $(I,g_I=\pm {\rm
d}t^2)$ is conformal Killing. Indeed, if $X$ is a vector field on
$(I,g_I)$, then $X$ can be expressed as $X = h \partial_t$ for
some smooth function $h \in \mathcal{C}^\infty(I)$.

For the rest of the paper, let $M=I_f \times F$ be a standard static
space-time with the metric $g=f^2 g_I \oplus g_F$, where $g_I=-{\rm d}t^2$.
Suppose that $X,Y,Z \in \mathfrak X(I)$ and $V,W,U \in \mathfrak
X(F)$, then (see \cite{U01})
\begin{equation} \label{eq:Lie deriv 1a}
{\rm L}_{X+V}g(Y+W,Z+U)  =  f^2 {\rm L}^I_X g_I(Y,Z) + 2f
V(f)g_I(Y,Z)
 +  {\rm L}^F_V g_F(W,U).
\end{equation}
Moreover, we also have
\begin{equation} \label{eq:Lie deriv 2}
{\rm L}_{h \partial_t} g_I(Y,Z) = Y(h)g_I(Z,\partial_t) + Z(h) g_I
(Y, \partial_t).
\end{equation}
By combining (\ref{eq:Lie deriv 1a}) and (\ref{eq:Lie deriv 2}),
we can state the following result.
\begin{theorem} \label{con-1} \cite{DU08a} Let $M=I_f \times F$ be a standard
static space-time with the metric $g=-f^2{\rm d}t^2 \oplus g_F.$
Suppose that $h \in C^{\infty}(I)$ and $V \in \mathfrak X(F)$.
Then $h
\partial_t + V$ is a conformal-Killing vector
field on $M$ with $\sigma \in C^{\infty}(M)$ if and only if the
following properties are satisfied:
\begin{enumerate}
\item $V$ is conformal-Killing on $F$ with associated $\sigma \in
C^{\infty}(F),$
\item $h$ is affine, i.e, there exist $\mu, \nu \in \mathbb{R}$
such that $h(t)=\mu t + \nu$ for any $t \in I,$
\item $V(f) = (\sigma-\mu) f.$
\end{enumerate}
\end{theorem}
\noindent Consequently, $h \partial_t + V$ is a Killing vector
field on $M$ if and only if the following properties are
satisfied:
\begin{enumerate}
\item $V$ is Killing on $F$,
\item there exist $\mu, \nu \in \mathbb{R}$
such that $h(t)=\mu t + \nu$ for any $t \in I,$
\item $V(f) = - \mu f$.
\end{enumerate}

In \cite{DU08a}, to provide a characterization of Killing vector
fields on standard static space-times, we modify the procedure
used in \cite{S99} (see also \cite{CC93}) to study the structure
of Killing and conformal-Killing vector fields on warped products.
In \cite{S99}, the author obtains full characterizations of the
Killing and conformal-Killing vector fields on generalized
Robertson-Walker space-times.
Here, we will state some of the main results about the characterization
of Killing vector fields obtained in \cite{DU08a}.

Let $(F,g_F)$ be a Riemannian manifold of dimension $s$ admitting
at least one \textit{nonzero} Killing vector field. Thus, there
exists a basis $\{K_{\overline b} \in \mathfrak X(F) | \,
\overline b = 1,\cdots, \overline s \}$ for the set of Killing
vector fields on $F$.
%
%\noindent
At this point, we would like to emphasize that the
dimension of the set of conformal Killing vector fields on
$(I,-{\rm d}t^2)$ is infinite, so that one cannot apply directly
the procedure in \cite{S99} before observing that the form of
conformal Killing vector fields on $(I,-{\rm d}t^2)$ is trivial
(i.e, any vector field on $(I,-{\rm d}t^2)$ is conformal). Adapting
the S\'{a}nchez technique to $M=I _f\times F$, a vector
field $K \in \mathfrak X(M)$ is a Killing vector field if and only
if $K$ can be written in the form
\begin{equation}\label{eq:Killing structure}
   K= \psi h \partial_t +
\phi^{\overline b} K_{\overline b},
\end{equation}
where $h$, $\phi^{\overline b} \in C^\infty(I)$ for any $\overline
b \in\{1, \cdots, \overline m\}$ and $\psi \in C^\infty(F)$
satisfy
\begin{equation} \label{2eqnst}
\left\{
\begin{array}{rcl}
 h^\prime \psi+ \phi^{\overline
b} K_{\overline b}(\ln f)& = & 0 \\
{\rm d} \phi^{\overline b} \otimes g_F(K_{\overline b},\cdot) +
g_I(h \partial_t, \cdot) \otimes f^2 {\rm d} \psi & = & 0.
\end{array}
\right.
\end{equation}
Since $\displaystyle{{\rm d} \phi^{\overline b} = (\phi^{\overline
b})^\prime {\rm d}t}$ with $\phi^{\overline b} \in \mathcal
C^\infty(I)$ and $g_I(h \partial_t, \cdot)=-h {\rm d}t$,
\eqref{2eqnst} is equivalent to
\begin{equation} \label{3eqnst}
\left\{
\begin{array}{rcl}
 h^\prime \psi+ \phi^{\overline
b} K_{\overline b}(\ln f)& = & 0 \\
(\phi^{\overline b})^\prime {\rm d}t \otimes g_F(K_{\overline
b},\cdot) & = & h {\rm d}t \otimes f^2 {\rm d} \psi .
\end{array}
\right.
\end{equation}

\noindent The following notation will be useful. Let $h$ be a
continuous function defined on a real interval $I$. If there
exists a point $t_0 \in I$ such that $h(t_0)\neq 0$, then
$I_{t_0}$ denotes the connected component of $\{t \in I: h(t)\neq
0\}$ such that $t_0 \in I_{t_0}$.

\noindent By the method of separation of variables and a detailed
analysis of system \eqref{3eqnst}, one can state the following result.
\begin{theorem}\label{thm:Killing ssst} \cite{DU08a} Let $(F,g_F)$
be a Riemannian manifold, $f \in C^\infty_{>0}(F)$ and
$\{K_{\overline b}\}_{1 \le \overline b \le \overline m}$ a basis of
Killing vector fields on $(F,g_F)$. Let also $I$ be an open interval
of the form $I=(t_1,t_2)$ in $\mathbb R,$ where
$-\infty \leq t_1 < t_2 \leq \infty.$
Consider the standard static space-time $I _f\times F$ with the
metric $g=-f^2{\rm d}t^2 \oplus g_F$.

\noindent Then, any Killing vector field on $I _f\times F$ admits
the structure
\begin{equation}\label{eq:eq:Killing 2-s no 0}
K = \psi h \partial_t + \phi^{\overline b} K_{\overline b}
\end{equation}
where $h$ and $\phi^{\overline b} \in C^\infty(I)$ for any
$\overline b \in\{1, \cdots, \overline m\}$ and $\psi \in
C^\infty(F).$

\noindent Furthermore, assume that $K$ is a vector field on $I
_f\times F$ with the structure as in \eqref{eq:eq:Killing 2-s no
0}. Hence,
\begin{itemize}
\item [{\bf (i)}] if $h \equiv 0,$ then the vector field
$K=\phi^{\overline b} K_{\overline b}$ is Killing on  $I _f\times
F$ if and only if the functions $\phi^{\overline b}$ are constant
and $\phi^{\overline b} K_{\overline b}(\ln f)=0$.
\item [{\bf (ii)}] if $h \equiv h_0 \neq 0$ is constant,
then the vector field $K= \psi h_0 \partial_t + \phi^{\overline b}
K_{\overline b}$ is Killing on  $I_f \times F$ if and only if is
satisfied
\begin{equation}
\label{3eqnst vector field h=h0 2} \left\{
\begin{array}{l}
%f \in C^\infty_{>0}(F),\psi \in C^\infty(F);\\
 f^{2} {\rm grad}_F \psi \textrm{ is a Killing vector field on } (F,g_F)
 \textrm{ with }\\
\textrm{ coefficients } \{\tau_{\overline b}\}_{1 \le {\overline
b} \le {\overline m}} \textrm{ relative to the basis }
\{K_{\overline b}\}_{1 \le
{\overline b} \le {\overline m}};\\
(f^{2} {\rm grad}_F \psi)(\ln f)=0 \,(i.e, \,{\rm grad}_F \psi(f)=0); \\
\forall {\overline b}: \phi^{\overline b}(t)=h_0 \tau^{\overline
b}t+\omega^{\overline b} \textrm{ with } \omega^{\overline b}\in
\mathbb{R} : \omega^{\overline b}
K_{\overline b}(\ln f)=0.\\
\end{array}
\right.
\end{equation}
\item [{\bf (iii)}] if
%\begin{equation}\label{eq:eq:Killing 2-s no 0}
%K = \psi h \partial_t + \phi^{\overline b} K_{\overline b}
%\end{equation}
$K$ is a Killing vector field on $I _f\times F$ with the
nonconstant function $h$, then the set of functions $h$, $\psi$
and $\{\phi_{\overline b}\}_{1 \le \overline b \le \overline m}$
satisfy
\begin{equation}\label{eq:3eqnst 2 h no 0}
\left\{
\begin{array}{l}
(a)\left\{
\begin{array}{l}
%f \in C^\infty_{>0}(F),
\psi \equiv 0; \\
\phi^{\overline b}(t)=\omega^{\overline b} \textrm{ on } I_{t_0}
\textrm{ where }\omega^{\overline b} \in
\mathbb{R}: \omega^{\overline b} K_{\overline b}(\ln f)=0\\
\end{array}
\right.
\\
\textrm{or }\\
(b)\left\{
\begin{array}{l}
%f \in C^\infty_{>0}(F),\psi \in C^\infty(F);\\
 f^{2} {\rm grad}_F \psi \textrm{ is a Killing vector field on }
(F,g_F)\\
\textrm{with } \textrm{coefficients } \{\tau_{\overline b}\}_{1
\le {\overline b} \le {\overline m}} \, \textrm{ relative to the } \\
\textrm{basis } \{K_{\overline b}\}_{1 \le
{\overline b} \le {\overline m}};\\
(f^{2} {\rm grad}_F \psi)(\ln f)  =  \nu \psi \textrm{ where }
\nu \, \textrm{is constant} ;\\
    h(t)= \left\{
    \begin{array}{lcl}
    a e^{\sqrt{-\nu}\,t} + b e^{-\sqrt{-\nu}\,t}\, \,
    & \textrm{ if } & \nu \neq 0 \\
    a t + b \, \,  &\textrm{ if }&  \nu = 0,\\
    \end{array}
    \right.\\
    \textrm{with } a, b \in \mathbb{R};\\
\forall {\overline b}: \phi^{\overline b}(t)=\tau^{\overline
b}\displaystyle \int_{t_0}^t h(s)ds +\omega^{\overline b} \textrm{
with }
\omega^{\overline b}\in \mathbb{R}:\\
\displaystyle h^\prime (t_0) \psi + \omega^{\overline
b}K_{\overline b}(\ln f)=0 \textrm{ on }
I_{t_0}\\
\end{array}
\right.
\\
\end{array}
\right.
\end{equation}
for any $t_0 \in I$ with $h(t_0)\neq 0.$

\noindent Conversely, if a set of functions $h$, $\psi$ and
$\{\phi_{\overline b}\}_{1 \le \overline b \le \overline m}$,
satisfy \eqref{eq:3eqnst 2 h no 0} with an arbitrary $t_0$ in $I$
and the entire interval $I$ (instead of $I_{t_0}$) and $\psi \in
C^\infty(F)$, then the vector field $\tilde{K}$ on the standard
static space-time $I _f\times F$ associated to the set of
functions as in \eqref{eq:eq:Killing 2-s no 0} is Killing on $I
_f\times F$.
\end{itemize}
\end{theorem}

For clarity, we also state the following lemma which covers the
case where the Riemannian manifold $(F,g_F)$ admits no nonidentical
zero Killing vector field.

\begin{lemma}\label{lem:0 Killing}
Let $(F,g_F)$ be a Riemannian manifold of dimension $s$ and $f \in
C^\infty_{>0}(F)$. Let also $I$ be an open interval of the form
$I=(t_1,t_2)$ in $\mathbb R,$ where $-\infty \leq t_1 < t_2 \leq
\infty$. Suppose that the only Killing vector field on $(F,g_F)$ is
the zero vector field. Then all the Killing vector fields on the
standard static space-time $I_f\times F$ are given by
$h_0 \partial_t$ where $h_0$ is a constant.
\end{lemma}

\textit{Theorem \ref{thm:Killing ssst}} is relevant to the problem
given by:

\begin{equation}
\label{eq:pb Killing 1} \left\{
\begin{array}{l}
f \in C^\infty_{>0}(F), \psi \in C^\infty(F);\\
f^{2} {\rm grad}_F \psi \textrm{ is a Killing vector field on }
(F,g_F);\\
(f^{2} {\rm grad}_F \psi)(\ln f)  =  \nu \psi, \nu \in \mathbb{R} .\\
\end{array}
\right.
\end{equation}
%and define $\mathcal{K}_f^\nu := \{ \psi \in C^\infty(F): \psi
%\textrm{ verifies } \eqref{eq:pb Killing 1} \}$.
We are interested in the existence of nontrivial solutions for
\eqref{eq:pb Killing 1}. To study this, for any $Z \in \mathfrak
{X}(F)$ and $\varphi \in C^\infty(F)$ we define the (0,2)-tensor
on $(F,g_F)$ given by
\begin{equation}\label{eq:special Killing 2}
    B_Z^\varphi(\cdot,\cdot):=\textrm{d} \varphi (\cdot)
    \otimes g_F(Z,\cdot) + g_F(\cdot,Z) \otimes \textrm{d}
    \varphi (\cdot).
\end{equation}

A central role in our study of \eqref{eq:pb Killing 1} is played
by the next proposition which also shows up the relevance of the
Hessian tensor in all these questions.

\begin{prp} \label{prp:special Killing} \cite{DU08a}
Let $(F,g_F)$ be a Riemannian manifold, $f \in C^\infty_{>0}(F)$
and $\psi \in C^\infty(F)$. Then the vector field $f^{2} {\rm
grad}_F \psi$ is Killing on $(F,g_F)$ if and only if
\begin{equation}\label{eq:special Killing 1}
   {\rm H}_F^\psi + \frac{1}{f} B_{{\rm grad}_F \psi}^{f}=0.
\end{equation}
\end{prp}
The latter proposition and the identity $fg_F({\rm grad}_F
\psi,{\rm grad}_F f)=(f {\rm grad}_F \psi)(f)$, allow to express
\eqref{eq:pb Killing 1} in the equivalent form
\begin{equation} \label{eq:3eqnst 4}
\left\{
\begin{array}{l}
f \in C^\infty_{>0}(F), \psi \in C^\infty(F);\\
\displaystyle {\rm H}_F^\psi + \frac{1}{f}
B_{{\rm grad}_F \psi}^{f} = 0 ;\\
fg_F({\rm grad}_F \psi,{\rm grad}_F f) = \nu \psi, \nu \in
\mathbb{R}.
\end{array}
\right.
\end{equation}

%\begin{remark} \label{rem:Lie algebra 1}
By \textit{Proposition \ref{prp:special Killing}}, if the
dimension of the Lie algebra of Killing vector fields of $(F,g_F)$
is zero, then the system \eqref{eq:3eqnst 4} has only the trivial
solution given by a constant $\psi $ (this constant is not $0$
only if $\nu =0$). This happens, for instance when $(F,g_F)$ is a
compact Riemannian manifold of negative-definite Ricci curvature
without boundary, indeed it is sufficient to apply the vanishing
theorem due to Bochner (see for instance \cite{B46}, \cite[Theorem
1.84]{B}).
% or \cite[Proposition 6.6 of Chapter III]{ST}).
%\end{remark}

The next \textit{Lemma \ref{lem: laplace spectrum 1}} allows to
prove  that the system \eqref{eq:3eqnst 4} is still
\textit{equivalent} to
\begin{equation} \label{eq:3eqnst 5}
\left\{
\begin{array}{l}
f \in C^\infty_{>0}(F), \psi \in C^\infty(F);\\
\displaystyle {\rm H}_F^\psi + \frac{1}{f}
B_{{\rm grad}_F \psi}^{f} = 0 ;\\
-\Delta_{g_F} \psi = \displaystyle \nu \frac{2}{f^2}
    \psi \textrm{ where } \nu
\textrm{ is a constant}.
\end{array}
\right.
\end{equation}

\begin{lemma}\label{lem: laplace spectrum 1}\cite{DU08a}
Let $(F,g_F)$ be a Riemannian manifold and $f \in
C^\infty_{>0}(F)$. If $(\nu,\psi)$ satisfies \eqref{eq:3eqnst 4},
then $\nu$ is an eigenvalue and $\psi$ is an associated
$\nu-$eigenfunction of the elliptic problem:
\begin{equation}\label{eq:weight Laplace-Beltrami 1}
    -\Delta_{g_F} \psi = \nu \frac{2}{f^2}
    \psi \, \textrm{ on } (F,g_F).
\end{equation}
\end{lemma}

Thus, by arguments of critical points and maximum principle, we
obtain the following characterization results.

\begin{prp}\label{prp:laplace spectrum 1} Let $(F,g_F)$ be a compact
Riemannian manifold and $f \in C^\infty_{>0}(F)$. Then $(\nu,\psi)$
satisfies \eqref{eq:3eqnst 4} if and only if $\nu =0$ and $\psi$ is
constant.
\end{prp}

\begin{theorem}\label{thm:killing compact fiber} Let $M=I _f\times F$
be a standard static space-time with the metric
$g=-f^2{\rm d}t^2 \oplus g_F.$ If $(F,g_F)$ is compact then, the set
of all Killing vector fields on the standard static space-time $(M,g)$
is given by
$$\{a
\partial_t + \tilde{K}| \, a \in \mathbb{R}, \tilde{K} \textrm{ is a
Killing vector field on } (F,g_F) \textrm{ and } \tilde{K}(f)=0
\}.$$
\end{theorem}

\begin{example} \label{rem:sharipov2007}
\textit{(Killing vector fields in the Einstein static universe)}
In \cite{Sh07}, the author studied Killing vector fields of a
closed homogeneous and isotropic universe (for related questions
in quantum field theory and cosmology see \cite{F87, LL}). Theorem
6.1 of \cite{Sh07} corresponds to Theorem \ref{thm:killing
compact fiber} for the spherical universe $\mathbb{R} \times
\mathbb{S}^3$ with the pseudo-metric $\displaystyle -( R^2 {\rm
d}t^2 -R^2 h_0)$, where the sphere $\mathbb{S}^3$ endowed with the
usual metric $h_0$ induced by the canonical Euclidean metric of
$\mathbb{R}^4$ and $R$ is a real constant (i.e., a stable
universe).
\end{example}

As we have already mentioned, any Killing vector field of a
compact Riemannian manifold of negative-definite Ricci tensor is
equal to zero. Thus, one can easily state the following result.
\begin{cor} \label{cor:killing compact fiber of neggative definite
Ricci tensor} Let $M=I _f\times F$ be a standard static space-time
with the metric $g=-f^2{\rm d}t^2 \oplus g_F.$ Suppose that
$(F,g_F)$ is a compact Riemannian manifold of negative-definite
Ricci tensor. Then, any Killing vector field on the standard static
space-time $(M,g)$ is given by $a \partial_t$ where $a \in \mathbb
R.$
\end{cor}

In \cite[Theorem 5]{S07}, it is shown that the decomposition of a
space-time as a standard static one is essentially unique when the
fiber $F$ is compact. We observe that {\it Corollary
\ref{cor:killing compact fiber of neggative definite Ricci
tensor}} enables us to establish a stronger conclusion (i.e.,
nonexistence of a nontrivial strictly stationary
\footnote{Here, a stationary field means that it is Killing and time-like
at the same time (see \cite{S07}).}
field) under a stronger assumption involving the definiteness of
the Ricci tensor.

At this point, we would like to make some comments about the case
where the Riemannian part of a standard static space-time is not
compact. While the Theorem \ref{thm:Killing ssst} does not require
the compactness of the Riemannian manifold $(F,g_F)$, this assumption
is the central idea for a complete characterization similar
to the one in Theorem \ref{thm:killing compact fiber}. The key question
in our approach is the full characterization of the solutions of
\eqref{eq:3eqnst 5} (or the equivalent problems \eqref{eq:pb Killing 1}
and \eqref{eq:3eqnst 4}) which is reached if $(F,g_F)$ is compact. In the
noncompact case, the latter question is more difficult. It is possible to
obtain partial nonexistence results for \eqref{eq:3eqnst 5}, but
the global question is still open. However, there are particular
situations, like Example \ref{rem:sharipov2007}, where the
application of Theorem \ref{thm:Killing ssst} is sufficient for a
complete classification.

Other relevant and related problem is the full classification of
the conformal Killing vector fields of a standard static
space-time. There are partial recent results in this direction
(see for instance \cite{AC, Sha-Iq} and the references therein).

\section{Energy Conditions}

Recall that a space-time is said to satisfy the \emph{strong
energy condition}, briefly SEC, if ${\rm Ric}(X,X) \geq 0$ for all
causal tangent vectors $X$ and the \emph{time-like}(respectively,
\emph{null
%, space-like
}) \emph{convergence condition}, briefly TCC
(respectively, NCC
%, SCC
), if ${\rm Ric}(X,X) \geq 0$ for all
time-like (respectively, null
%, space-like
) tangent vectors $X$.
Notice that the SEC implies the NCC. Furthermore the TCC is
equivalent to the SEC, by continuity.
%However, we want to
%emphasize that they are not equivalent in terms of physical
%implications. Indeed, the
The actual difference between TCC and SEC follows from the fact
that while TCC is just a geometric condition imposed on the Ricci
tensor, SEC is a condition on the stress-energy tensor. They can
be considered equivalent due to the Einstein equation (see below
\eqref{eq: Einstein eq}).

Moreover, a space-time is called to satisfy the \emph{weak energy
condition}, briefly WEC, if ${\rm T}(X,X) \geq 0$ for all
time-like vectors, where ${\rm T}$ is the energy-momentum tensor,
which is determinated by physical considerations.

\noindent Along this article, when we consider the energy-momentum
tensor, we will assume that the Einstein equation holds
(see \cite{HE,ON}). More explicitly,
\begin{equation}\label{eq: Einstein eq} {\rm Ric}-\frac{1}{2} \tau g
= 8 \pi {\rm T}.
\end{equation}

%Notice that in particular \eqref{eq: Einstein eq} gives the
%explicit form of the energy-momentum tensor $T$.

\noindent The WEC has many applications in general relativity
theory such as nonexistence of closed time-like curve (see
\cite{ChoPark}) and the problem of causality violation
(\cite{OriSoen}). But its fundamental usage still lies in
Penrose's Singularity theorem (see \cite{Pen1}).

Let $M=I _f\times F$ be a standard static space-time with the
metric $g=-f^2{\rm d}t^2 \oplus g_F$. By \eqref{eq:Qricci-sss},
${\rm Ric} (\partial_t,\partial_t)= f \Delta_F f$. So, since
$g(\partial_t,\partial_t)= -f^2 <0$, the warping function $f$ is
necessarily subharmonic, i.e. $\Delta_F f \geq 0$, if the standard
static space-time satisfies the SEC (or equivalently, the
TCC)\cite{ DA1}. On the other hand, it is well known that there is
no nonconstant subharmonic functions on compact Riemannian
manifolds \cite{B46}, and hence $f$ is a positive constant if $(F,g_F)$
is compact. Furthermore, applying a family of Liouville type
results of Li, Schoen and Yau, in \cite{DU08b} we give a set of
sufficient conditions implying the warped function is a positive
constant under the hypothesis that $(F,g_F)$ is complete and
noncompact.

Below, we state a set of necessary conditions for a standard
static space-time to satisfy the NCC and other Ricci curvature
conditions which are useful to study conformal hyperbolicity
through the studies of Markowitz, more precisely Theorems 5.1 and
5.8 in \cite{MM1}. We also observe that there are more accurately
analogous results for a Generalized Robertson-Walker space-time
given in \cite[Proposition 4.2]{ES00} (see also \cite{DU08b}).

\begin{theorem} \label{ec-t} \cite{DU08b,DA1} Let $M=I _f\times F$
be a standard static space-time with the metric $g=-f^2{\rm d}t^2 \oplus g_F$,
where $s = \dim F \ge 2$.
\begin{enumerate}
\item If $\mathcal{R}ic_F$ and $\mathcal{Q}^f_F$ are positive
semi-definite, then $M$ satisfies the TCC and the NCC.

\item If $\mathcal{R}ic_F$ and $\mathcal{Q}^f_F$ are negative
semi-definite, then ${\rm Ric}(w,w) \leq 0$ for any causal
vector $w \in \mathfrak X(M)$.

\item If $(F,g_F)$ is Ricci flat, then $\mathcal{Q}^f_F$ is positive
semi-definite if and only if $M$ satisfies the NCC.
\end{enumerate}
\end{theorem}

Now we state a small results about energy conditions in terms of
the energy-momentum tensor $T$. It is easy to obtain from
\eqref{eq: Einstein eq}, \eqref{eq:Qricci-sss} and
\eqref{eq:sc-sss} that for any $U \in \mathfrak X (I)$ and $V \in
\mathfrak X (F)$ is
\begin{equation}\label{eq: em tensor}
    8\pi T(U+V,U+V)= \mathcal{R}ic_F(V) + \frac{1}{f}
    \mathcal{Q}_F^f(V) - \frac{1}{2} \tau_F g(U+V,U+V).
\end{equation}
So, as above, since $g(\partial_t,\partial_t)= -f^2 <0$ results
that if a standard static space-time satisfies the WEC, then
$\tau_F \geq 0$ and as consequence
\begin{theorem} \label{ec-w} \cite{DU08b} Let $M=I _f\times F$ be
a standard static space-time with the metric $g=-f^2{\rm d}t^2 \oplus g_F$,
where $s = \dim F \ge 2$.
\begin{enumerate}
\item If $\mathcal{R}ic_F$ and $\mathcal{Q}^f_F$ are positive
(respectively negative) semi-definite, then ${\rm T}(w,w) \geq 0$
(respectively $\leq 0$) for any causal vector $w \in \mathfrak X(M).$

\item If $(F,g_F)$ is Ricci flat, then for any $u \in \mathfrak X (I)$
and $v \in \mathfrak X (F)$ is $8\pi
T(u+v,u+v)=\mathcal{Q}^f_F(v)$. Thus, $\mathcal{Q}^f_F$is positive
semi-definite if and only if ${\rm T}(w,w) \geq 0$ for any vector
$w \in \mathfrak X(M)$.
\end{enumerate}
\end{theorem}

%\subsection{Conformal Hyperbolicity}

In \cite{MM1} the intrinsic Lorentzian pseudo-distance
$d_M \colon M \times M \to [0,\infty)$ was defined by
\begin{equation}\label{eq:pseudo-distance}
d_M(p,q)=\inf_{\alpha}L(\alpha),
\end{equation}
where the infimum is taken over all the chains of null geodesic
segments joining $p$ and $q$ and $L(\alpha)$ means the length of
the chain $\alpha$. Such a chain $\alpha$ is a sequence of points
$p=p_0,p_1,\dots,p_k=q$ in $M,$ pairs of points $(a_1,b_1), \dots,
(a_k,b_k)$ in $(-1,1)$ and projective maps (i.e., a projective map
is simply a null geodesics with the projective parameter as the
natural parameter) $f_1,\dots,$ $f_k$ from $(-1,1)$ into $M$ such
that $f_i(a_i)=p_{i-1}$ and $f_i(b_i)=p_i$ for $i=1,\cdots,k$.
Besides, the length of $\alpha $ is
\begin{equation*}
    L(\alpha)=\sum_{i=1}^{k} \rho(a_i,b_i),
\end{equation*}
where $\rho$ is the Poincar\'{e} distance in $(-1,1)$ (see
\cite{MM1, MM2} for details). Notice that $d_M$ is really a
pseudo-distance, i.e., it is non-negative, symmetric and satisfies
the triangle inequality. A Lorentzian manifold $(M,g)$ where $d_M$
is a distance is called \textit{conformally hyperbolic}.

\noindent In \cite[Theorem 5.1]{MM1}, it is proved that if $(M,g)$
is a null geodesically complete Lorentz manifold satisfying the
reverse NCC condition, i.e. ${\rm Ric}(X,X) \leq 0$ for all null
vectors $X$, then it has a trivial Lorentzian pseudo-distance,
i.e., $d_M \equiv 0$. Moreover, in \cite[Theorem 5.8]{MM1}, it is
obtained that if $(M,g)$ is an $n(\geq 3)$-dimensional Lorentzian
manifold satisfying the NCC and the null generic condition,
briefly NGC, (i.e., ${\rm Ric}(\gamma^\prime,\gamma^\prime) \neq
0,$ for at least one point of each inextendible null geodesic
$\gamma$) then, it is conformally hyperbolic.

\noindent Under the light of these theorems, one can easily
conclude that
\begin{itemize}
\item Complete Einstein space-times (in particular, Minkowski,
de-Sitter and the anti-de Sitter space-times) have all trivial
Lorentzian pseudo-distances because of Theorem 5.1 of \cite{MM1}.

\item The Einstein static universe has also trivial Lorentzian
pseudo-distance since the space-times in the previous item can be
conformally imbedded in the Einstein static universe.

\item A Robertson-Walker space-time (i.e., an isotropic
homogeneous space-time) is conformally hyperbolic due to Theorem
5.9 of \cite{MM1}.

\item The Einstein-de Sitter space $M$ is conformally hyperbolic
and (see Theorem 5 in \cite{MM2} for details and a precise formula
for the Lorentzian pseudo-distance on this class of space-time).
\end{itemize}

Applying \textit{Theorem \ref{ec-t}}, \eqref{eq:Qricci-sss} and
the previous Markowitz results, we obtain the theorems that
follow.

%Now, we will state the following results which can be deduced from
%Theorem 5.1 in \cite{MM1}, \textit{Theorem \ref{ec-t}} and the
%null geodesic completeness of the underlying standard static
%space-times (see \cite[Theorem 3.12]{DA2}).

\begin{theorem} \label{main-12} \cite{DU08b}
Let $M=\mathbb R _f\times F$ be a standard static space-time with
the metric $g=-f^2{\rm d}t^2 \oplus g_F.$ Suppose that
$\mathcal{R}ic_F$ and $\mathcal{Q}^f_F$ are negative
semi-definite.
\begin{enumerate}
\item If $(F,g_F)$ is compact, then the Lorentzian pseudo-distance
$d_M$ on the standard static space-time $(M,g)$ is trivial, i.e.,
$d_M \equiv 0.$ \item If $(F,g_F)$ is complete and  $0 < \inf f $,
then the Lorentzian pseudo-distance $d_M$ on the standard static
space-time $(M,g)$ is trivial, i.e., $d_M \equiv 0.$
\end{enumerate}
\end{theorem}

In the Theorem \ref{main-12}, the general hypothesis ensure the
\textit{reversed} NCC, in order to apply \cite[\textsc{Theorem
5.1}]{MM1}. The additional hypothesis in item (2) of the same
theorem implies the null geodesic completeness of $M$ by
\cite[Thoerem 3.12]{ADt}. We observe that there are more general
hypotheses which imply null geodesic completeness, see for
instance \cite[Th. 3.9(ii b)]{RS}.

%We state the following result about conformal hyperbolicity
%of standard static space-times by combining Theorem 5.8 in
%\cite{MM1} and \textit{Theorem \ref{ec-t}}, and
%\eqref{eq: ricci vector}.

\begin{theorem} \label{main-3} \cite{DU08b}
Let $M=I _f\times F$ be a standard static space-time with the
metric $g=-f^2{\rm d}t^2 \oplus g_F$. Suppose that
$\mathcal{R}ic_F$ is positive semi-definite and $\mathcal{Q}^f_F$
is positive definite. Then the standard static space-time $(M,g)$
is conformally hyperbolic.
\end{theorem}

Now we state some results joining the conformal hyperbolicity and
causal conjugate points of a standard static space-time by using
\cite{BJ1,BE1,BE2,CE} and also \cite{BEE}. In \cite[Theorem
2.3]{CE}, it was shown that if the line integral of the Ricci
tensor along a complete causal geodesic in a Lorentzian manifold
is positive, then the complete causal geodesic contains a pair of
conjugate points.

Assume that $\gamma=(\alpha, \beta)$ is a complete causal geodesic
in a standard static space-time of the form $M=I _f\times F$ with
the metric $g=-f^2{\rm d}t^2 \oplus g_F.$ Then by using
$g(\gamma^\prime, \gamma^\prime) \leq 0$ and \eqref{eq:Qricci-sss}
we have,
\begin{equation*}\label{} {\rm Ric}(\gamma^\prime,
\gamma^\prime)= {\rm Ric}_F(\beta^\prime, \beta^\prime) +
\frac{1}{f}Q^f_F(\beta^\prime,\beta^\prime)-
\underbrace{g(\gamma^\prime,\gamma^\prime)}_{\leq 0}
\frac{1}{f}\Delta_F f.
\end{equation*}

We can easily state the following existence result for conjugate
points of complete causal geodesics in a conformally hyperbolic
standard static space-time by {\it Theorem \ref{main-3}} and
\cite[Theorem 2.3]{CE}.

\begin{theorem} \label{cor-1} Let $M=I _f\times F$ be a standard
static space-time with the metric $g=-f^2{\rm d}t^2 \oplus g_F$.
Suppose that $\mathcal{R}ic_F$ is positive semi-definite. If
$\mathcal{Q}^f_F$ is positive definite, then $(M,g)$ is
conformally hyperbolic and any complete causal geodesic in $(M,g)$
has a pair of conjugate points.
\end{theorem}

By using \eqref{eq:Qricci-sss}, \cite[Propositions 11.7, 11.8 and
Theorem 11.9]{BEE} and \cite[Corollary 3.17]{ADt}, we can
establish an existence result for conjugate points of time-like
geodesics in a standard static space-time which by \textit{Theorem
\ref{cor-1}} is also conformally hyperbolic.
%%%%%
%Here, the Lorentzian arc-length of a time-like geodesic $\gamma$
%is denoted by ${\rm L}(\gamma)$ (see page 135 of \cite{BEE}).
%Furthermore, the time-like diameter of a space-time $(M,g)$ is
%defined by
%$$ {\rm diam}(M,g) = \sup \{ d(p,q) | \, p,q \in M \}$$ where
%$d(p,q)$ is the Lorentzian distance function given in terms
%of the Lorentzian arc-length (see pages 137 and 399 of
%\cite{BEE} for further details).
%%%%%

\noindent In the next theorem, $\mathbf{L}$ denotes the usual
time-like Lorentzian length and ${\rm diam}_\mathbf{L}$ denotes
the corresponding time-like diameter (see \cite[Chapters 4 and
11]{BEE}).

\begin{theorem} \label{cor-3} Let $M=I _f\times F$ be a standard
static space-time with the metric $g=-f^2{\rm d}t^2 \oplus g_F$.
Suppose that $\mathcal{R}ic_F$ and $\mathcal{Q}^f_F$ are positive
semi-definite. If there exists a constant $c$ such that
$\displaystyle\frac{1}{f}\Delta_F f \geq c >0$, then
\begin{enumerate}
\item any time-like geodesic $\gamma \colon [r_1,r_2] \to M$ in $(M,g)$
with $\mathbf{L}(\gamma) \geq \pi \sqrt{\frac{n-1}{c}}$ has a pair
of conjugate points,
\item for any time-like geodesic $\gamma \colon [r_1,r_2] \to M$ in
 $(M,g)$ with $\mathbf{L}(\gamma)
> \pi\sqrt{\frac{n-1}{c}}$,  $r=r_1$ is conjugate along $\gamma$ to some
$r_0 \in (r_1,r_2),$ and as consequence $\gamma$ is not maximal,
\item if $I=\mathbb{R}$, $(F,g_F)$ is complete and $\sup f<\infty$,
 then $\displaystyle{\rm diam}_\mathbf{L}(M,g) \leq
 \pi\sqrt{\frac{n-1}{c}}.$
\end{enumerate}
\end{theorem}

%Note that, by using Theorem 7.1 of \cite{MM1}, one can also deduce
%that the group of conformal automorphisms of the the underlying
%standard static space-time has a compact isotropy group at each
%point $p$ when the hypothesis in {\it Theorem \ref{main-3}}, {\it
%Theorems \ref{cor-1}, \ref{cor-2} or \ref{cor-3}} are verified.

In the final part of \cite{DU08b} we show some examples and
results connecting the tensor $Q_F^f$, conformal hyperbolicity,
concircular scalar fields and Hessian manifolds, where the role of
the Hessian tensor is central.

\begin{center}
\textsc{Acknowledgements}
\end{center}
The authors wish to thank the referee for the useful and
constructive suggestions. F. D. thanks The Abdus Salam
International Centre of Theoretical Physics for their warm
hospitality where part of this work has been done.

\bibliographystyle{amsalpha}

\end{document}